\documentclass[11pt,reqno]{amsart}
\usepackage{fullpage,amsfonts,amsmath,amscd,amssymb}

\theoremstyle{plain}

\newtheorem*{lem}{Lemma}
\newtheorem*{prop}{Proposition}
\newtheorem*{thm}{Theorem}

\newtheorem*{cor}{Corollary}

\theoremstyle{remark}

\newcommand{\ep}{\epsilon}

\DeclareMathOperator{\hil}{Hilb^n(\mathbb{C}^2)}
\DeclareMathOperator{\symn}{Sym^n(\mathbb{C}^2)}

\DeclareMathOperator{\mat}{Mat} 
\DeclareMathOperator{\edo}{End} 
 \DeclareMathOperator{\spc}{Spec}
\DeclareMathOperator{\homo}{Hom} \DeclareMathOperator{\im}{im}

 \DeclareMathOperator{\md}{mod}
\DeclareMathOperator{\Ind}{Ind}

\title{Baby Verma modules for rational Cherednik algebras}
\author{Iain Gordon}
\email{ig@maths.gla.ac.uk}
\thanks{}
\begin{document}
\begin{abstract}
In this paper, we introduce \textit{baby Verma modules} for
symplectic reflection algebras of complex reflection groups at
parameter $t=0$ (the so--called rational Cherednik algebras at
parameter $t=0$) and present their most basic properties.  As an
example, we use baby Verma modules to answer several problems
posed by Etingof and Ginzburg, \cite{EG}, and give an elementary
proof of a theorem of Finkelberg and Ginzburg, \cite{FG}.
\end{abstract}
\maketitle
\section{Introduction}
\subsection{}
Symplectic reflection algebras algebras of complex reflection
groups (the so--called rational Cherednik algebras) arise in many
different mathematical disciplines: integrable systems, Lie
theory, representation theory, differential operators, symplectic
geometry. Early on a dichotomy (depending on a single parameter
$t$) appears in the behaviour of these algebras: when $t\neq 0$
the algebras are rather non--commutative and relations with
differential operators come to the fore; when $t=0$ the algebras
have very large centres and connections between the geometry of
the centre and representation theory are of interest. This paper
concentrates on the second case. We introduce a family of finite
dimensional modules called \textit{baby Verma modules} and present
their most basic properties. By analogy with the representation
theory of reductive Lie algebras in positive characteristic, we
believe these modules are fundamental to the understanding of the
representation theory and associated geometry of the rational
Cherednik algebras at parameter $t=0$.
\subsection{}\label{res1} Let $W$ be a complex reflection group and $\mathfrak{h}$ its reflection
representation over $\mathbb{C}$. The rational Cherednik algebras
introduced by Etingof and Ginzburg, \cite{EG}, are deformations of
the skew--group ring $\mathbb{C}[\mathfrak{h}\oplus
\mathfrak{h}^*]\ast W$ depending on parameters $t$ and $c$. We
denote them by $H_{t,c}$. It can be shown that when $t=0$ there is
an embedding of $\mathbb{C}[\mathfrak{h}^*]^W \otimes
\mathbb{C}[\mathfrak{h}]^W$ into $Z_c$, the centre of $H_{0,c}$,
and thus a map
$$\Upsilon : \spc (Z_c) \longrightarrow \mathfrak{h}^*/W \times
\mathfrak{h}/W.$$ The baby Verma modules are $H_{0,c}$--modules
naturally associated to points in the (reduced) zero fibre,
$\Upsilon^{-1}(0)$. Studying their properties allows us to deduce
some properties of $Z_c$ and of the map $\Upsilon$. We prove the
following results:
\begin{itemize} \item There is a surjective map, $\Theta$, from the isomorphism classes of irreducible representations
of $W$ over $\mathbb{C}$ to the elements in the fibre
$\Upsilon^{-1}(0)$ .
\item The map $\Theta$ is an isomorphism if and only if
the points of $\Upsilon^{-1}(0)$ are smooth in $\spc (Z_c)$.
\item If $\Theta$ is an isomorphism and $S$ a simple $W$--module, then the component of the scheme theoretic
fibre $\Upsilon^*(0)$ corresponding to $\Theta(S)$ has a
$\mathbb{C}^*$--action whose associated Poincar\'{e} polynomial is
explicitly described. The component has, in particular, dimension
$\dim (S)^2$.
\item The variety $\spc (Z_c)$ is singular for $W$ a finite Coxeter group of
type $D_{2n}, E, F, H$ and $I_2(m)$ ($n\geq 2, m\geq 5$).
\end{itemize}
These results are related to problems of Etingof and Ginzburg,
confirming \cite[Conjecture 17.14]{EG}, answering \cite[Question
17.15]{EG} and partially answering \cite[Question 17.1]{EG}.
\subsection{}
The final result above should be compared to recent results on the
non--existence of a crepant resolution for the orbit space
$(\mathfrak{h}^*\times \mathfrak{h})/W$. It is shown in \cite{kal}
that the orbit space admits no crepant resolution if $W$ is the
Coxeter group of type $G_2$; folding arguments then imply no such
resolution for $W$ of type $D$ and $E$. It is an interesting
problem to make the relationship between resolutions and
deformations of $(\mathfrak{h}\times\mathfrak{h}^*)/W$ precise.
\subsection{} Recall Calogero--Moser
phase space defined by $$\mathcal{CM}_n= \{ (X,Y)\in
\mat_n(\mathbb{C})\times \mat_n(\mathbb{C}) : [X,Y] + Id =
\text{rank 1 matrix}\}//PGL_n(\mathbb{C}). $$ In the particular
case of $W=\mathfrak{S}_n$, the symmetric group on $n$ letters, it
is known that for $c\neq 0$ there is an isomorphism between $\spc
(Z_c)$ and $\mathcal{CM}_n$ under which the map $\Upsilon$
corresponds to sending the pair of matrices $(X,Y)$ to their
eigenvalues. Since it is known that $\mathcal{CM}_n$ is smooth,
our results mentioned in \ref{res1} associate to each point in
$\Upsilon^{-1}(0)$ a partition of $n$, as in \cite{wil}. Let
$\lambda\vdash n$. We prove:
\begin{itemize}
\item The Poincar\'{e} polynomial of the component of the scheme--theoretic
fibre $\Upsilon^*(0)$ corresponding to $\lambda$ is
$K_{\lambda}(t)K_{\lambda}(t^{-1})$, where $K_{\lambda}$ is the
Kostka polynomial corresponding to $\lambda$.
\item For the point $M_{\lambda} \in \Upsilon^{-1}(0)$
corresponding to $\lambda$, there is a relationship between the
finite dimensional algebra $H_{0,c}/M_{\lambda}H_{0,c}$ and
Springer theory for the nilpotent orbit in $\mat_n(\mathbb{C})$
with Jordan normal form of type $\lambda$.
\end{itemize}
These results confirm a pair of conjectures of Etingof and
Ginzburg, \cite[Conjectures 17.12 and 17.13]{EG}, and give an
elementary proof of a recent theorem of Finkelberg and Ginzburg,
\cite{FG}. As noted later, the extension of the first result to
$\Gamma\wr \mathfrak{S}_n$, the wreath product of a finite
subgroup of $SL_2(\mathbb{C})$ by the symmetric group, is
straightforward.
\subsection{} The paper is organised as follows. Sections 2--3
introduce notation and present the embedding of
$\mathbb{C}[\mathfrak{h}^*]^W \otimes \mathbb{C}[\mathfrak{h}]^W$
into $Z_c$. In Section 4 we define and study baby Verma modules.
Section 5 sees the proofs of all but the last result stated in
\ref{res1}, whilst Section 6 deals with the special case of the
symmetric group, $\mathfrak{S}_n$. In Section 7 we prove the final
result of \ref{res1} and present a couple of basic open problems.
\subsection{Acknowledgements} The author is grateful to Alexander Kuznetsov
 for explaining to me the relationship between fixed
points in Calogero--Moser spaces and Hilbert schemes used at the
end of Section \ref{kuz}, and to Ken Brown for allowing me to
include his unpublished result in \ref{brown}. The author is
partially supported by the Nuffield Foundation grant NAL/00625/G.
\section{Notation}
\subsection{} Let $W$ be a complex reflection group and
$\mathfrak{h}$ its reflection representation over $\mathbb{C}$.
Let $\mathcal{S}$ denote the set of complex reflections in $W$.
Let $\omega$ be the canonical symplectic form on $V =
\mathfrak{h}\oplus \mathfrak{h}^*$. For $s\in\mathcal{S}$, let
$\omega_s$ be the skew--symmetric form which coincides with
$\omega$ on $\im (\text{id}_V-s)$ and has $\ker (\text{id}_V-s)$
as the radical. Let $c: \mathcal{S} \longrightarrow \mathbb{C}$ be
a $W$--invariant function sending $s$ to $c_s$. The rational
Cherednik algebra at parameter $t=0$ (depending on $c$) is the
quotient of the skew group algebra of the tensor algebra on $V$,
$TV \ast \mathbb{C}W$, by the relations
\[
[x,y ] = \sum_{s\in \mathcal{S}} c_s \omega_s(x,y) s
\]
for all $y\in \mathfrak{h}^*$ and $x\in\mathfrak{h}$. This algebra
is denoted by $H_c$.
\subsection{} We will denote a set of
representatives of the isomorphism classes of simple $W$--modules
by $\Lambda$.
\subsection{} Let $R=\oplus R_i$ be a finite dimensional
$\mathbb{Z}$-graded $\mathbb{C}$--algebra. We will let $R-\md$
($R-\md_{\mathbb{Z}}$) denote the category of finite dimensional
(graded) left $R$--modules. Given an object $M$ in
$R-\md_{\mathbb{Z}}$, let $M[i]$ denote the $i$th shift of $M$,
defined by $M[i]_j = M_{j-i}$. The \textit{graded endomorphism
ring} of $M$ is defined as $$\edo(M) = \bigoplus_{i\in \mathbb{Z}}
\homo_{R-\md_{\mathbb{Z}}}(M, M[i]).$$ Given $N$, a
$\mathbb{Z}$--graded right $R$--module, the \textit{graded dual}
$N^{\circledast}$ is the object of $R-\md_{\mathbb{Z}}$ defined by
$$ (N^{\circledast})_i = \{ f: N \longrightarrow \mathbb{C} : f(
N_j) = 0 \text{ for all }j\neq -i \}.$$ The forgetful functor is
denoted $F: R-\md_{\mathbb{Z}}\longrightarrow R-\md$.
\section{The centre of $H_c$}
\subsection{}
Throughout $Z_c = Z(H_c)$ and
$A=\mathbb{C}[\mathfrak{h}^*]^W\otimes
\mathbb{C}[\mathfrak{h}]^W$. Thanks to \cite[Proposition
4.15]{EG}, $A\subset Z_c$ for any finite complex reflection group
$W$ and any parameter $c$. We will give an elementary proof of
this inclusion.
\subsection{}
Let $s\in \mathcal{S}$. There exists $\alpha_s\in \mathfrak{h}$
such that $\mathfrak{h} = \mathfrak{h}_s \oplus
\mathbb{C}\alpha_s$, where $\mathfrak{h}_s= \{ h\in\mathfrak{h}:
s(h)=h\}$ and $s(\alpha_s) = e^{2\pi in_s}\alpha_s$ for some
rational number $0 < n_s < 1$. Similarly, we have
$\alpha_s^{\vee}$, $-n_s$ and $\mathfrak{h}_s^*$. We don't choose
to normalise $\alpha_s$ and $\alpha_s^{\vee}$, but note that for
all $h\in \mathfrak{h}_s$ it follows from
$$ \alpha_s^{\vee}( h ) = \alpha_s^{\vee} (s(h)) = e^{2\pi
in_s}\alpha_s^{\vee}(h)$$ that $\alpha_s^{\vee}(\mathfrak{h}_s) =
0$. Similarly, $\alpha_s (\mathfrak{h}^*_s) = 0$.
\subsection{}
 For any
reflection $s\in W$, an explicit calculation shows, for
$y\in\mathfrak{h}^*$ and $x\in\mathfrak{h}$
$$ \omega_s (x,y) =
\frac{\alpha_s(y)\alpha_s^{\vee}(x)}{\alpha_s^{\vee}(\alpha_s)}.$$
The principal commutation relation for $H_c$ is
\begin{equation}\label{cheredrel} [x,y] = \sum_{s\in S}
c_s\frac{\alpha_s(y)\alpha_s^{\vee}(x)}{\alpha_s^{\vee}(\alpha_s)}s.\end{equation}

\subsection{} Let $s$ be a reflection in $W$ and set $\lambda_s =1-e^{2\pi in_s}\in\mathbb{C}$.
For $h\in \mathfrak{h}$ $$ s(h) = h - \lambda_s
\frac{\alpha_s^{\vee}(h)}{\alpha_s^{\vee}(\alpha_s)} \alpha_s.$$
This action extends multiplicatively to
$\mathbb{C}[\mathfrak{h}^*]$.
\subsection{} The functional
$\alpha_{s}^{\vee}$ extends inductively to a $\mathbb{C}$--linear
operator on $\mathbb{C}[\mathfrak{h}^*]$ as follows:
$\alpha_s^{\vee}|_{\mathbb{C}} \equiv 0$ and for $p,p'\in
\mathbb{C}[\mathfrak{h}^*]$, \begin{equation}\label{newop}
\alpha_{s}^{\vee}(pp') = \alpha_s^{\vee}(p)p' +
p\alpha_s^{\vee}(p') -
\lambda_s\frac{\alpha_s^{\vee}(p)\alpha_s^{\vee}(p')}{\alpha_s^{\vee}(\alpha_s)}
\alpha_s.\end{equation} Define a graded operator on
$\mathbb{C}[\mathfrak{h}^*]$ by \begin{equation}\label{tils}
\tilde{s}(p) = p - \lambda_s
\frac{\alpha_s^{\vee}(p)}{\alpha_s^{\vee}(\alpha_s)} \alpha_s.
\end{equation}
An inductive calculation proves that
\begin{equation}\label{sameops}\text{the operators $s$ and $\tilde{s}$ agree
on $\mathbb{C}[\mathfrak{h}^*]$}.\end{equation} In particular,
\begin{equation}\label{fixedpts}p\in \mathbb{C}[\mathfrak{h}^*]^W
\text{ if and only if } \alpha_s^{\vee}(p) = 0 \text{ for all
reflections }s\in W.\end{equation}
\subsection{} Let $y\in\mathfrak{h}^*$ and $p\in \mathbb{C}[\mathfrak{h}^*]$.
One proves inductively using (\ref{cheredrel}), (\ref{newop}),
(\ref{tils}) and (\ref{sameops}) that
$$ [p,y] = \sum_{s} c_s \frac{y(\alpha_s)\alpha_s^{\vee}(p)}
{\alpha_s^{\vee}(\alpha_s)} s. $$ It follows from (\ref{fixedpts})
that if $p\in\mathbb{C}[\mathfrak{h}^*]^W$ then $p\in Z_c$. Of
course similar arguments apply to $\mathbb{C}[\mathfrak{h}]^W$, so
we have proved
\begin{prop}
\label{ctre} Let $H_c$ be a Cherednik algebra. Then $A\subset
Z_c$.
\end{prop}
\noindent We denote the corresponding morphism of varieties
\begin{equation}\label{ups} \Upsilon :\spc (Z_c) \longrightarrow \spc
(A)= \mathfrak{h}^*/W \times \mathfrak{h}/W.
\end{equation}
\section{Baby Verma modules}
\subsection{}\label{defh} Recall by Proposition \ref{ctre} that
$A\subset Z_c$. Let $A_+$ denote the elements of $A$ of with no
scalar term. Define
$$\overline{H_c} = \frac{H_c}{A_+H_c}.$$
This is a $\mathbb{Z}$--graded algebra where $\deg(x) = 1$
(respectively $\deg(y)=-1$, $\deg(w)= 0$) for $x\in \mathfrak{h}$
(respectively $y\in\mathfrak{h}^*, w\in W$). The PBW theorem,
\cite[Theorem 1.3]{EG}, shows that $\overline{H_c}$ has a (vector
space) triangular decomposition \begin{equation}
\label{tri}\begin{CD} \overline{H_c} @> \sim
>> \mathbb{C}[\mathfrak{h}^*]^{co W} \otimes
\mathbb{C}[\mathfrak{h}]^{co W}  \otimes \mathbb{C}W,
\end{CD}\end{equation} where ``$co W$" stands for $W$--coinvariants, that is
$\mathbb{C}[\mathfrak{h}^*]^{co W} = \mathbb{C}
[\mathfrak{h}^*]/\mathbb{C}
[\mathfrak{h}^*]_+^W\mathbb{C}[\mathfrak{h}^*]$. Thanks to
Chevalley's Theorem these coinvariant rings are graded versions of
the regular representation of $W$, so in particular
$\dim(\overline{H_c}) = |W|^3$.
\subsection{}
Let $\overline{H_c^-} = \mathbb{C}[\mathfrak{h}]^{co W}\ast
\mathbb{C}W$, a subalgebra of $\overline{H_c}$. The algebra map
$\overline{H_c^-} \longrightarrow \mathbb{C}W$ sending an element
$q\ast w$ to $q(0)w$, makes any $W$--module into an
$\overline{H_c^-}$--module.

Given $S\in \Lambda$, a simple $W$--module, we define $M(S)$, a
\textit{baby Verma module}, by $$M(S) = \overline{H_c}
\otimes_{\overline{H_c^-}} S.$$ By construction $M(S)$ is an
object of $\overline{H_c}-\md_{\mathbb{Z}}$.  Equation (\ref{tri})
shows that $\dim(M(S)) = |W|\dim(S)$. Furthermore for $T\in
\Lambda$, set
\begin{equation} \label{fake} f_T(t) =
\sum_{i\in\mathbb{Z}}(\mathbb{C}[\mathfrak{h}^*]^{co W}:
T[i])t^i.\end{equation} The polynomials $f_T(t)$ are called the
\textit{fake degrees} of $W$ and have the property $f_T(1)=\dim
T$. They have been calculated for all finite Coxeter groups. In
the graded Grothendieck group of $W$, we have by construction
\begin{equation} \label{grt} [M(S)] = \sum_{T\in\Lambda}
f_T(t)[T\otimes S].\end{equation}
\subsection{} The following proposition gathers results
from \cite[Section 3]{NP}. Recall $F:
\overline{H_c}-\md_{\mathbb{Z}}\longrightarrow \overline{H_c}-\md$
is the forgetful functor.
\begin{prop}
\label{ind} Let $S,T\in \Lambda$.
\begin{enumerate}
\item The baby Verma module $M(S)$ has a
 simple head, denoted $L(S)$.
\item $M(S)$ is isomorphic to $M(T)[i]$ if and only if $S$ and
$T$ are the same element of $\Lambda$ and $i=0$.
\item $\{ L(S)[i] : S\in \Lambda, i\in\mathbb{Z}\}$ gives a complete set
of pairwise non--isomorphic simple graded
$\overline{H_c}$--modules.
\item $F(L(S))$ is a simple $\overline{H_c}$--module, and $\{
F(L(S)) : S\in\Lambda\}$ is a complete set of pairwise
non--isomorphic simple $\overline{H_c}$--modules.
\item Let $P(S)$ be the projective cover of $L(S)$. Then $F(P(S))$
is the projective cover $F(L(S))$.
\end{enumerate}
\end{prop}

\subsection{} \label{vermmul} \textbf{Lemma.} \textit{
Let $S, T\in \Lambda$ and $i$ be a positive integer. Then
$(M(S)[i]:L(T)) =0$ and $(M(S):L(T)) = \delta_{S,T}.$}
\begin{proof} By construction the baby Verma module $M(T)$
 is concentrated in non--negative degree and generated by its degree zero component,
$1\otimes T = M(T)_0$. Since $M(T)$ is a homomorphic image of
$P(T)$, we can find $\tilde{T}$, a $W$--submodule of $P(T)_0$
isomorphic to $T$ which maps onto $M(T)_0$. Since $P(T)$ has a
simple head, it follows that $\tilde{T}$ generates $P(T)$.

If $L(T)$ is a composition factor of $M(S)[i]$, there exists a
non--zero graded homomorphism from $P(T)$ to $M(S)[i]$, sending
$P(T)_0$ to $M(S)_{-i}$. As $M(S)_{-i}$ is zero for positive
values of $i$, this proves the first claim of the lemma. If $i$ is
zero, then $\tilde{T}$ maps to $M(S)_0$ and the second claim
follows.
\end{proof}

\subsection{} An object $V$ in $\overline{H_c}-\md_{\mathbb{Z}}$ is said to have an
\textit{$M$-filtration} if it has a filtration $$ 0=V_0 \subseteq
V_1 \subseteq \cdots \subseteq V_{n-1}\subseteq V_n=V$$ such that
for $1\leq j \leq n$, $V_j/V_{j-1} \cong M(S)[i]$ for $S$ a simple
$W$--module and $i\in \mathbb{Z}$. By \cite[Corollary 4.3]{NP} the
numbers $$ [V : M(S)[i] ] = |\{ j: V_j/V_{j-1}\cong M(S)[i]\}|$$
are independent of the filtration used.

\subsection{} Let $\overline{H_c^+} = \mathbb{C}[\mathfrak{h^*}]^{co
W}\ast \mathbb{C}W$. Given a simple $W$--module, we define
$$M^-(S) = (S^* \otimes_{\overline{H_c^+}}
\overline{H_c})^{\circledast}.$$

It is shown in \cite[Theorem 4.5]{NP} that for a simple
$W$--module $S$, the projective cover $P(S)$ has an $M$-filtration
such that for any $i\in\mathbb{Z}$ and simple $W$--module $T$
\begin{equation} \label{BGG1} [P(S): M(T)[i]] = (M^{-}(T)[i] : L(S)). \end{equation}
\subsection{} \label{bg} For this paragraph only
suppose that $W$ is a finite Coxeter group. Then (\ref{BGG1}) can
be improved to the following Brauer--type reciprocity formula
\begin{equation} \label{BGG2} [P(S): M(T)[i]] = (M(T)[i]: L(S)).\end{equation}
 This follows from an
adaptation of \cite[Theorem 5.1]{NP} with two special ingredients.
Firstly, we require the antiautomorphism, $\omega$, of $H_c$ which
sends $x\in \mathfrak{h}$, respectively $y\in\mathfrak{h}^*$, to
the corresponding element $\tilde{x}$ of $\mathfrak{h}^*$,
respectively $\tilde{y}$ of $\mathfrak{h}$, under the
$W$--invariant form $( -,-)$ on $\mathfrak{h}$ and $w$ to $w^{-1}$
for $w\in W$. Secondly, we must observe that any simple
$W$--module is self--dual (since the characters of $W$ all take
values in $\mathbb{R}$).
\section{On the fibre $\Upsilon^*(0)$}
\subsection{} Recall the morphism $\Upsilon$ given in (\ref{ups}). We
will study the scheme theoretic fibres $\Upsilon^*(0)$ of the
point $0\in \mathfrak{h}^*/W\times \mathfrak{h}/W$. The closed
points of these fibres, denoted $\Upsilon^{-1}(0)$, are precisely
the maximal ideals of $Z_c$ lying over $A_+$. Given such an ideal,
$M$, the component of $\Upsilon^*(0)$ corresponding to $M$ is
defined to be $\spc (\mathcal{O}_M)$, where
$$\mathcal{O}_M = (Z_c)_M/A_+(Z_c)_M .$$
The algebra $\mathcal{O}_M$ inherits a $\mathbb{Z}$--grading from
$\overline{H_c}$, seen as follows. The $\mathbb{Z}$--grading on
$\overline{H_c}$ corresponds to an algebraic action of
$\mathbb{C}^*$ which preserves both $Z_c$ and $A$ and makes
$\Upsilon$ a $\mathbb{C}^*$-equivariant morphism. Since $0$ is the
unique fixed point of $\mathfrak{h}^*/W\times \mathfrak{h}/W$,
each fibre $\Upsilon^*(0)$ inherits a $\mathbb{C}^*$--action, in
other words a $\mathbb{Z}$--grading.
\subsection{}
\label{azsm} Recall that a point $M\in \spc (Z_c)$ is called
\textit{Azumaya} if and only if there is a unique simple
$H_c$--module annihilated by $M$ and its dimension equals $|W|$.
By \cite[Theorem 1.7]{EG} the Azumaya points are precisely the
smooth points of $\spc (Z_c)$ and any simple $H_c$--module lying
over an Azumaya point is isomorphic, as a $W$--module, to the
regular representation of $W$.
\subsection{}
The algebra $\overline{H_c}$ of \ref{defh} splits into a direct
sum of indecomposable algebras, the blocks, which by
\cite[Corollary 2.7]{BG} are in one--to--one correspondence with
the elements of $\Upsilon^{-1}(0)$
$$\overline{H_c} = \bigoplus_{M\in \Upsilon^{-1}(0)}
\mathcal{B}_M.$$ If $M$ is an Azumaya point of $\spc (Z_c)$ there
is an isomorphism
\begin{equation} \label{az} \mathcal{B}_M \cong \mat_{|W|}
(\mathcal{O}_M),\end{equation} by \cite[Proposition 2.2]{BG}.
\subsection{} For any simple $W$--module $S$, the baby Verma module
$M(S)$ is a $\mathbb{Z}$--graded indecomposable
$\overline{H_c}$--module by Proposition \ref{ind}(1), and so a
non--trivial module for a unique block. This gives a mapping
\begin{equation} \label{the} \Theta :\{ \text{isomorphism classes of simple
$W$--modules}\} \longrightarrow \{ \text{elements of
$\Upsilon^{-1}(0)$}\}.\end{equation} This mapping is surjective,
since each simple $\overline{H_c}$--module is a quotient of a baby
Verma by Proposition \ref{ind}(3).
\subsection{}
\label{poin} Suppose that $M\in\Upsilon^{-1}(0)$ is smooth in
$\spc (Z_c)$. It follows from \ref{azsm} and (\ref{az}) that
\begin{equation*} \label{az2} \mathcal{B}_M \cong \mat_{|W|}
(\mathcal{O}_M).\end{equation*} In particular, $\mathcal{B}_M$ has
a unique simple module, implying that $M=\Theta(S)$ for a unique
element $S$ of $\Lambda$. Moreover, we have an isomorphism
\begin{equation*}\label{grim} \begin{CD} \mathcal{O}_{M}
@> \sim >> \edo(P(S)). \end{CD}
\end{equation*} The map is graded since it is given by multiplication
by elements of $\mathcal{O}_{M}$. Thus we have a formula for
$p_S(t)\in\mathbb{Z}[t,t^{-1}]$, the Poincar\'{e} polynomial
 of $\mathcal{O}_{M}$,
\begin{eqnarray*}  p_S(t) &=& \sum_{i \in
\mathbb{Z}} \dim \left(\homo_{R-\md_{\mathbb{Z}}}(P(S),
P(S)[i])\right)t^i \\ &=& \sum_{i\in \mathbb{Z}} (P(S)[i]:
L(S))t^i.\end{eqnarray*} Since the block $\mathcal{B}_{M}$ has
only the simple modules $L(S)[i]$ ($i\in \mathbb{Z}$) it follows
from (\ref{BGG1}) that
\begin{eqnarray} \label{pp}
 p_S(t) & = & \sum_{i,j\in \mathbb{Z}} [P(S)[i]: M(S)[j]] (M(S)[j]: L(S)) t^i \nonumber \\
 & = & \sum_{i,j\in\mathbb{Z}} (M^-(S)[j-i]:L(S))
(M(S)[j]: L(S)) t^{i-j}t^j \nonumber \\ &= & \left(
\sum_{i\in\mathbb{Z}} (M^-(S)[i]:L(S)) t^{-i} \right)\left(
\sum_{i\in\mathbb{Z}}(M(S)[i] : L(S))t^i \right).
\end{eqnarray}
\subsection{} Recall the fake polynomials, $f_S(t)$, defined in (\ref{fake}).
Let $b_S$ be the lowest power of $t$ appearing in $f_S(t)$.
 \begin{thm} \label{et} Suppose $M\in \Upsilon^{-1}(0)$ is smooth in $\spc
(Z_c)$. Then there exists a unique simple module $S\in \Lambda$
such that $\Theta (S) = M$. Furthermore $p_S(t)$, the Poincar\'{e}
polynomial of $\mathcal{O}_M$, is given by
$$p_S (t) = t^{b_{S^{\ast}}-b_S}f_S(t)f_{S^{\ast}}(t^{-1}).$$
\end{thm}
\begin{proof} Following \ref{poin}, it remains to prove the
formula for the Poincar\'{e} polynomial $p_S(t)$. By \ref{azsm}
$L(S)$ is isomorphic to the regular representation of $W$ and
$L(S)$ is the only possible composition factor of $P(S)$. We begin
by determining the unique integer $l_S$ for which $(L(S)_{l_S}:
\mathbf{1})\neq 0$, where $\mathbf{1}$ denotes the trivial
$W$--module. Let $j$ be the smallest integer such that
$(M(S)_{j}:\mathbf{1})\neq 0$. Since $L(S)$ is a homomorphic image
of $M(S)$, $j \leq l_S$. If $j < l_S$ then necessarily
$(M(S)[l_S-j]: L(S)) \neq 0$, contradicting Lemma \ref{vermmul}.
Hence $j=l_S$.

If $T$ is a simple $W$--module, $\mathbf{1}$ is a summand of
$T\otimes S$ if and only if $T\cong S^{\ast}$. Hence, by
(\ref{grt}),
$$\sum_{i\in\mathbb{Z}} (M(S)_k : \mathbf{1}) t^k = f_{S^{\ast}}(t).$$
Thus $l_S$ equals $b_{S^\ast}$, the lowest power of $t$ appearing
in $f_{S^{\ast}}(t)$. Moreover, for each copy of $\mathbf{1}$ in
$M(S)_k$ there occurs a corresponding copy of $L(S)[k-b_{S^\ast}]$
in $M(S)$, yielding $$\sum_{i\in\mathbb{Z}}(M(S)[i] :L(S))t^i =
\sum_{i\in\mathbb{Z}} (M(S): L(S)[-i])(t^{-1})^i =
t^{b_{S^{\ast}}}f_{S^{\ast}}(t^{-1}).$$ A similar argument shows
that
$$\sum_{i\in\mathbb{Z}}(M^-(S)[i] :L(S))t^{-i} = t^{-b_{S}}
f_{S}(t).$$ The result follows.
\end{proof}

\subsection{} \label{easyrems} We remark that the statement of Theorem \ref{et}
simplifies if $W$ is a finite Coxeter group. As observed in
\ref{bg} every simple $W$--module is self--dual, so the
Poincar\'{e} polynomial above becomes
$$p_S(t) = f_S(t)f_S(t^{-1}).$$

\subsection{}
On specialising $t$ to $1$ and recalling that $f_S(1)=\dim S$, we
obtain the following corollary, which generalises
\cite[Proposition 4.16]{EG}, solves \cite[Problem 17.15]{EG} and
incidentally confirms \cite[Conjecture 17.14]{EG}.
\begin{cor}
Suppose $\Upsilon^{-1}(0)$ consists of smooth points in
$\spc(Z_c)$. Then $\Theta: \Lambda \longrightarrow
\Upsilon^{-1}(0)$ is a bijection such that the scheme theoretic
multiplicity of $\Upsilon^{-1}(0)$ at the point $\Theta (S)$ is
$\dim(S)^2$.
\end{cor}

\section{The symmetric group case}\label{kuz}
\subsection{} Throughout this section $W$ will be the symmetric group
$\mathfrak{S}_n$ and $H_{c}$ will be the Cherednik algebra of $W$
with \textit{non--zero} parameter $c$. This choice of $c$ ensures
that $\spc(Z_c)$ is smooth, \cite[Corollary 16.2]{EG}.

\subsection{} We recall notation and several basic facts on the representation
theory of $\mathfrak{S}_n$ over $\mathbb{C}$. The simple
$\mathfrak{S}_n$--modules are indexed by partitions of $n$, which
we will write as $\lambda = (\lambda_1,\lambda_2,\ldots )$. For
such a partition $\lambda \vdash n$, we denote the simple module
by $S_{\lambda}$. In particular $S_{(n)} = \mathbf{1}$ and
$S_{(1^n)} = \mathbf{sign}$, the trivial and sign representation
of $\mathfrak{S}_n$ respectively. Given $\lambda$, a partition of
$n$, we let $\lambda'\vdash n$ denote the transpose of $\lambda$.
We have $S_{\lambda'} \cong S_{\lambda}'$, where, for any
$\mathfrak{S}_n$--module $V$, we set $V'= V\otimes \mathbf{sign}$.
The Young subgroup $\mathfrak{S}_{\lambda} \leq \mathfrak{S}_n$ is
the row stabiliser of the Young tableau of shape $\lambda$, with
the tableau numbered in the order that one reads a book. As an
abstract group it is isomorphic to the product $\prod_i
\mathfrak{S}_{\lambda_i}$.

\subsection{} We recall two polynomials which arise in the theory of
symmetric functions. The (two variable) \textit{Kostka--Macdonald
coefficients} $K_{\lambda\mu}(q,t)$ are the transition functions
between the Macdonald polynomials and Schur functions, see
\cite[VI.8]{mac}. The \textit{Kostka polynomials} are defined as
$K_{\lambda\mu}(0,t)$. In particular, we have
$$K_{\lambda}(t) \equiv t^{-n(\lambda)}K_{\lambda(1^n)}(0,t)= (1-t)\cdots (1-t^n) \prod_{u\in \lambda}
(1-t^{h_{\lambda}(u)})^{-1} \in \mathbb{Z}[t,t^{-1}], $$ where
$n(\lambda') = \sum_{i\geq 1}(i-1)\lambda_i'$ and $h_{\lambda}(u)$
is the hook--length of $u$ in $\lambda$.
\subsection{}
The first part of the following theorem recovers \cite[Theorem
1.2]{FG}.
\begin{thm} \label{maincor} Let $\lambda$ be a partition of $n$. \\ $(1)$ The Poincar\'{e}
 polynomial of $\mathcal{O}_{\Theta(S_{\lambda})}$ is given by
$$p_{S_{\lambda}}(t) = K_{\lambda}(t)K_{\lambda}(t^{-1}).$$
\noindent $(2)$ The image of the simple $H_c$--module
$L(S_{\lambda})$ in the graded Grothendieck group of $W$ is
$$[L(S_\lambda)] = \sum_{\mu\vdash n}
K_{\mu\lambda}(t,t)[S_\mu].$$
\end{thm}
\begin{proof} Let $\mu, \rho$ be partitions of $n$. We let $H_{\mu}(t) =
 \prod_{u\in\mu} (1-t^{h_{\mu}(u)})$ denote the hook length polynomial. We will write $\chi^{\lambda}$
 for the character of $\mathfrak{S}_n$
 corresponding to $S_{\lambda}$, and $\chi^{\lambda}_{\rho}$ for the evaluation of $\chi^{\lambda}$ on
 a conjugacy class with cycle type $\rho$.

\noindent (1) By \cite[Theorem 3.2]{ste}
\begin{equation}\label{gch} f_{S_{\lambda}}(t) = (1-t)\cdots (1-t^n)
s_{\lambda}(1,t,t^2,\ldots ),
\end{equation} where $s_{\lambda}$ is the Schur function for
partition $\lambda$. The result follows from \ref{easyrems} and
\cite[Example I.3.2]{mac} which shows that
\begin{equation}\label{schur} s_{\lambda}(1,t,t^2,\ldots ) = t^{n(\lambda)}
H_{\lambda}(t)^{-1}.\end{equation}

\noindent (2) Let $p_{\lambda,\mu}(t) = \sum_i [L(S_{\lambda})_i:
S_\mu]t^i.$ In the Grothendieck group of graded $W$--modules we
have, by the proof of Theorem \ref{et},
\begin{eqnarray*}
[M(S_{\lambda})] &=&\sum_{\mu}
\left( \sum_i [M(S_{\lambda})_i : S_\mu]t^i\right) [S_{\mu}] \\
&=& \sum_{\mu} \left( \sum_{i,j} [M(S_{\lambda})[j]:
L(S_\lambda)][L(S_{\lambda})_{i+j}: S_{\mu}] t^i\right) [S_{\mu}]
\\ &=&  \sum_{\mu} t^{-b_{\lambda}}f_{S_\lambda}(t)
p_{\lambda,\mu}(t)[S_{\mu}].\end{eqnarray*} By Lemma
\ref{vermmul}, the transition matrix between the baby Verma
modules and the simple $\overline{H_c}$--modules is invertible in
$\mathbb{C}[[t]]$, so it is enough to show that \begin{equation}
\label{simmult} [M(S_{\lambda})] = \sum_{\mu}
t^{-b_{\lambda}}f_{S_\lambda}(t) K_{\mu \lambda}(t,t)
[S_{\mu}].\end{equation} We will prove (\ref{simmult}) by showing
that the characters of both sides are equal.

By (\ref{grt}) the character of $[M(S_{\lambda})]$ is $\sum_{\mu}
f_{S_{\mu}}(t)\chi^\mu \chi^\lambda,$ so we need to show that for
all partitions $\rho$
$$ \sum_{\mu}f_{S_\mu}(t)\chi^\mu_{\rho}\chi^\lambda_\rho=\sum_{\mu}
t^{-b_{\lambda}}f_{S_\lambda}(t)K_{\mu\lambda}(t,t)\chi^{\mu}_{\rho}.$$
By \cite[p.355]{mac}
\begin{equation}\label{import} \sum_{\mu} K_{\mu\lambda}(t,t)
\chi^{\mu}_{\rho} = \frac{H_{\lambda}(t)}{\prod_{i\geq 1}
(1-t^{\rho_i})}\chi^{\lambda}_{\rho},
\end{equation} and by (\ref{schur}) $$f_{S_\mu}(t) =
t^{n(\mu)}\prod_{j=1}^n
(1-t^j)H_\mu(t)^{-1}.$$ Note in particular this shows that $b_\mu
= n(\mu)$ for all partitions $\mu$. Therefore we have to prove
that
$$ \sum_{\mu}
f_{S_\mu}(t)\chi_\rho^\mu\chi_\rho^\lambda= \frac{\prod_{j=1}^n
(1-t^j)}{\prod_{i\geq 1}(1-t^{\rho_i})}\chi^\lambda_\rho $$ The
particular case of (\ref{import}) with $\lambda = (n)$ reduces the
above equation to \begin{equation} \label{fional} \sum_\mu
f_{S_\mu}(t)\chi_\rho^\mu \chi_\rho^\lambda= \sum_\mu
K_{\mu,(n)}(t,t) \chi_\rho^\mu\chi_\rho^{\lambda}.\end{equation}
By \cite[p.82 and p.362]{mac} $K_{\mu,(n)}(t,t)= f_{S_\mu}(t)$,
proving that indeed (\ref{fional}) holds.
\end{proof}
\subsection{} In \cite[Theorem 1.5]{FG} it is shown that an analogue of
Theorem
 \ref{maincor}(1)
 holds for the wreath product
$W=(\mathbb{Z}/N\mathbb{Z})\wr \mathfrak{S}_n$. This result can be
proved as above, using the fake degrees as described in
\cite[(5.5)]{ste}.
\subsection{Springer theory}\label{spring} Before proving Theorem \ref{3conj} we need to
 recall several results from Springer theory. Conjugacy
classes of nilpotent matrices in $\mat_n(\mathbb{C})$ are
classified by partitions of $n$ thanks to the Jordan normal form:
for each partition $\lambda\vdash n$, fix a representative of the
corresponding nilpotent conjugacy class, $e_{\lambda}$. Let
$\mathcal{B}$ be the flag variety of the Lie algebra
$\mat_n(\mathbb{C})$, consisting of all Borel subalgebras of
$\mat_n(\mathbb{C})$, and let $\mathcal{B}_{\lambda}\subseteq
\mathcal{B}$ denote the Borel subalgebras containing
$e_{\lambda}$. (Note that $\mathcal{B}= \mathcal{B}_{(1^n)}$.) Let
$d_{\lambda}= 2\dim_{\mathbb{C}}\mathcal{B}_{\lambda}$. Let
$H^*(\mathcal{B}_{\lambda},\mathbb{C})$ denote the singular
cohomology of $\mathcal{B}_{\lambda}$ with complex coefficients.
The following theorem presents the results we will use.
\begin{thm} For $\lambda\vdash n$ and $i\in \mathbb{N}$, there is an action of
$\mathfrak{S}_n$ on $H^i(\mathcal{B}_{\lambda},\mathbb{C})$
enjoying the following properties.
\begin{enumerate}
\item There is a $\mathfrak{S}_n$--equivariant isomorphism of algebras $\mathbb{C}[\mathfrak{h^*}]^{co \mathfrak{S}_n}
\cong H^*(\mathcal{B},\mathbb{C})$ doubling degree.
\item The simple module $S_{\lambda}$ is not a component of $H^i(\mathcal{B},\mathbb{C})$ for $i< d_{\lambda}$.
\item As an $\mathfrak{S}_n$--module, $H^{d_{\lambda}}(\mathcal{B}_{\lambda},\mathbb{C}) \cong S_{\lambda}$, whilst
 for odd $i$ or $i>d_{\lambda}$ $H^i(\mathcal{B}_\lambda,\mathbb{C}) =0$.
\item As an $\mathfrak{S}_n$--module,
$H^*(\mathcal{B}_{\lambda},\mathbb{C}) \cong
\Ind_{\mathfrak{S}_{\lambda}}^{\mathfrak{S}_n} \mathbf{1}.$
\item The inclusion $\mathcal{B}_{\lambda}\subseteq \mathcal{B}$
induces a surjective $\mathfrak{S}_n$--equivariant graded algebra
homomorphism
$$\pi_{\lambda}: H^*(\mathcal{B},\mathbb{C}) \longrightarrow
H^*(\mathcal{B}_{\lambda},\mathbb{C}).$$
\item Under the
$\mathbb{C}[\mathfrak{h}^*]^{co W}\ast \mathfrak{S}_n$--module
structure induced by $\pi_{\lambda}$,
$H^*(\mathcal{B}_{\lambda},\mathbb{C})$ has socle $S_{\lambda}$.
\end{enumerate}\end{thm}
\begin{proof} For (1)--(5) see
\cite{spr1},\cite{spr2},\cite{hotspr} and \cite{decpro}.

(6) Define a partial order on monomials in $\mathbb{C}[x_1,\ldots
,x_n]$ by declaring $x_1^{i_1}\ldots x_n^{i_n}\leq x_1^{j_1}\ldots
x_n^{j_n}$ if and only if $i_k\leq j_k$ for all $k=1,\ldots , n$.
By \cite[Proposition 4.2]{garpro}, the algebra
$\mathbb{C}[\mathfrak{h}^*]^{co W}/\ker \pi_\lambda$ has a basis
$\mathfrak{B}(\lambda)$ consisting of (the images of) monomials of
$\mathbb{C}[\mathfrak{h}^*] = \mathbb{C}[x_1,\ldots , x_n]$ with
the following property: the elements of $\mathfrak{B}(\lambda)$
are the monomials which are smaller in the above partial order
than the monomials in $\mathfrak{B}(\lambda)$ of degree
$d_\lambda$. We see therefore that the socle of
$H^{\ast}(\mathcal{B}_{\lambda},\mathbb{C})$ as a
$\mathbb{C}[\mathfrak{h}^*]^{co W}$--module is concentrated in
degree $d_{\lambda}$. The result follows from (3).
\end{proof}
In this situation, the \textit{Springer correspondence} is simply
the association of the nilpotent class parametrised by $\lambda$
to the simple $\mathfrak{S}_n$--module $S_{\lambda}\cong
H^{d_{\lambda}}(\mathcal{B}_{\lambda},\mathbb{C})$.
\subsection{}
Let $M_{\lambda} = \Theta (S_{\lambda}) \in \Upsilon^{-1}(0)$. By
\ref{azsm} we have an isomorphism
\begin{equation}\label{esim}\psi: \frac{H_{c}}{M_{\lambda}H_c}
\longrightarrow \edo_{\mathbb{C}}(L(S_\lambda))\end{equation}
defined by the action of $H_{c}$ on $L(S_{\lambda})$. Let
$\mathbf{e}, \mathbf{e}_-\in\mathbb{C}\mathfrak{S}_n\cong
L(S_\lambda)$ denote the trivial and sign idempotents
respectively. Fix $\epsilon_{\lambda}$ (respectively
$\epsilon_{\lambda}'$), a base vector in the one--dimensional
space $\mathbf{e} (H_c/M_{\lambda}H_c)\mathbf{e}_-$ (respectively
$\mathbf{e}_-(H_c/M_{\lambda}H_c)\mathbf{e}$). The subalgebras
$\mathbb{C}[\mathfrak{h}^*], \mathbb{C}[\mathfrak{h}] \subset H_c$
generate two subspaces $\mathbb{C}[\mathfrak{h}^*]\cdot
\epsilon_{\lambda}, \epsilon_{\lambda}\cdot
\mathbb{C}[\mathfrak{h}] \subset H_c/M_{\lambda}H_c$ which are
invariant under the left, respectively right
$\mathfrak{S}_n$-action on $H_c/M_{\lambda}H_c$ by multiplication.

The following theorem confirms \cite[Conjecture 17.12]{EG}.
\begin{thm}\label{3conj} In the above notation there are $\mathfrak{S}_n$--module isomorphisms $$
\mathbb{C}[\mathfrak{h}^*]\cdot \epsilon_{\lambda} \cong
\Ind_{\mathfrak{S}_{\lambda'}}^{\mathfrak{S}_n} \mathbf{1} \qquad
, \qquad \epsilon_{\lambda}\cdot \mathbb{C}[\mathfrak{h}] \cong
\Ind_{\mathfrak{S}_\lambda}^{\mathfrak{S}_n} \mathbf{sign}.$$ In
particular, these two $\mathfrak{S}_n$--modules have a single
non--zero irreducible component in common: it is isomorphic to
$S_{\lambda'}$.
\end{thm}
\begin{proof}
Throughout the proof let $R=\mathbb{C}[\mathfrak{h}^*]^{co
\mathfrak{S}_n}$. Let $N=1/2n(n-1)=d_{(1^n)}$. By Theorem
\ref{spring} $R_j=0$ if $j>N$ and $R_N$ is the sign
representation. Let $0\neq s\in R_N$. Recall that given $V$, a
$\mathfrak{S}_n$--module, we set $V' = V\otimes \mathbf{sign}$.

Under the isomorphism $\psi$ of (\ref{esim}), the elements
$\mathbf{e}$ and $\mathbf{e}_-$ correspond to projection onto the
trivial and sign components of $L(S_\lambda)$ respectively. Thus
$\epsilon_{\lambda}$ (respectively $\epsilon_{\lambda}'$) is the
unique (up to scalars) transformation which sends the sign
(respectively trivial) component of $L(S_{\lambda})$ to the
trivial (respectively sign) component of $L(S_{\lambda})$ and is
zero on all other isotypic components. As a left
$\mathfrak{S}_n$--module therefore,
$\mathbb{C}[\mathfrak{h}^*]\cdot \epsilon_{\lambda}$ (respectively
$\mathbb{C}[\mathfrak{h}^*]\cdot \epsilon_{\lambda}'$) corresponds
to the subspace $X=R\cdot t_{\lambda} \subseteq L(S_{\lambda})$
(respectively $Y=R \cdot s_{\lambda}$) where $t_{\lambda}$
(respectively $s_{\lambda}$) is a base vector for the trivial
(respectively sign) representation in $L(S_{\lambda})$.

By \cite[Theorem 1.7(ii)]{EG}, $H_c$ and $Z_c$ are Morita
equivalent, so in particular their injective dimensions both equal
$\dim \spc (Z_c)= 2(n-1)$. Since $\overline{H}_c$ is obtained from
$H_c$ by factoring out a regular sequence in $Z_c$ of length
$2(n-1)$ it follows that the injective dimension of
$\overline{H}_c$ is $0$. It follows that the projective cover
$P(S_{\lambda})$, and thus $M(S_{\lambda})$, has socle
$L(S_{\lambda})$. Therefore $t_{\lambda}$ (respectively
$s_{\lambda}$) belongs to the tensor product $S_{\lambda}\otimes
S_{\lambda} \subset L(S_{\lambda})\subseteq M(S_{\lambda})$
(respectively $S_{\lambda'}\otimes S_{\lambda}\subset
L(S_{\lambda}) \subseteq M(S_\lambda)$) where the first tensorand
is the highest degree component of $R$ isomorphic to $S_{\lambda}$
(respectively $S_{\lambda'}$). To be explicit, let $\{ v_i \}$ be
a basis of $S_{\lambda}$ and let $\{ f_i\}$ (respectively
$\{g_i\}$) be the corresponding basis for the copy of
$S_{\lambda}$ (respectively $S_{\lambda'}$) in the highest
possible degree of $R$. We take
$$t_{\lambda} = \sum_i f_i \otimes v_i \qquad (\text{respectively } s_{\lambda} = \sum_i g_i \otimes v_i).$$

Consider the map
$$\tau: R \longrightarrow X=R\cdot t_{\lambda}, \qquad (\text{respectively } \tau'
:R' \longrightarrow Y=r\cdot s_{\lambda},)$$ sending $p$ to
$p\cdot t_{\lambda}$ (respectively $p\otimes 1$ to $p\cdot
s_{\lambda}$). Since $t_{\lambda}\in \mathbf{1}$ (respectively
$s_{\lambda}\in \mathbf{sign}$), $\tau$ (respectively $\tau'$) is
a $\mathfrak{S}_n$--equivariant surjective homomorphism.

We prove first that $\tau$ (respectively $\tau'$) factors through
$\pi_{\lambda'}: H^*(\mathcal{B},\mathbb{C}) \longrightarrow
H^*(\mathcal{B}_{\lambda'},\mathbb{C})$ (respectively
$\pi_{\lambda}\otimes \text{id}: H^*(\mathcal{B},\mathbb{C})'
\longrightarrow H^*(\mathcal{B}_{\lambda},\mathbb{C})'$). Recall
the non--degenerate form $$( - , -):R \times R \longrightarrow
\mathbb{C},$$ defined by setting $(p,q)$ equal to the coefficient
of $pq$ in $R_N$ with respect to $s$. The form pairs $R_i$ with
$R_{N-i}$ and, since $(wr,ws)= \textbf{sign}(w)(r,s)$, a
representation $V$ with $V'$. Let $p\in \ker \pi_{\lambda'}$
(respectively $q\in \ker \pi_{\lambda}$). Observe that since
$\pi_{\lambda'}(Rp)=0$ (respectively $\pi_{\lambda}(Rq) = 0$) the
intersection of $Rp$ (respectively $Rq$) with the
$S_{\lambda'}$--component of
$H^{d_{\lambda'}}(\mathcal{B},\mathbb{C})$ (respectively the
$S_{\lambda}$--component of
$H^{d_{\lambda}}(\mathcal{B},\mathbb{C})$) is zero, by Parts (3)
and (5) of Theorem \ref{spring}. Thus $(R, pf_i) = (Rp, f_i) = 0$
(respectively $(R, qg_i) = (Rq, g_i) = 0$) for all $i$ . This
confirms that $\tau$ (respectively $\tau'$) factors through
$\pi_{\lambda'}$ (respectively $\pi_{\lambda}$).

We have an induced $R\ast \mathfrak{S}_n$--homomorphism
$${\sigma} : H^*(\mathcal{B}_{\lambda'},\mathbb{C}) \longrightarrow X, \qquad
(\text{respectively } {\sigma'} :
H^*(\mathcal{B}_{\lambda},\mathbb{C})' \longrightarrow Y).$$

Thanks to the non--degeneracy of the form $(-, -)$, we can find
$f_i'\in R$ such that $f_i'f_j =\delta_{ij}s\in R_N$. Then, for
any $i$, $f_i'\cdot t_{\lambda} = s \otimes s_i \in X$. Taking the
span over all $i$ shows that $S_{\lambda'} \cong R_N\otimes
S_\lambda \subseteq X$. It follows from Parts (2) and (6) of
Theorem \ref{spring} that the socle of
$H^*(\mathcal{B}_{\lambda'},\mathbb{C})$ does not lie in the
kernel of $\sigma$, and so $\sigma$ is injective. By Part (4) of
Theorem \ref{spring} we deduce that $$
\Ind_{\mathfrak{S}_{\lambda'}}^{\mathfrak{S}_n} \mathbf{1} \cong
H^*(\mathcal{B}_{\lambda'},\mathbb{C}) \cong X \cong
\mathbb{C}[\mathfrak{h}^*]\cdot \epsilon_{\lambda}.$$

Arguing as in the above paragraph we see that
$$\Ind_{\mathfrak{S}_{\lambda}}^{\mathfrak{S}_n} \mathbf{sign}
\cong H^*(\mathcal{B}_{\lambda},\mathbb{C})' \cong Y \cong
\mathbb{C}[\mathfrak{h}^*]\cdot \epsilon_{\lambda}'.$$ Let
$\omega$ be the antiautomorphism of $H_{c}$ discussed in \ref{bg}.
Since $\omega(M_{\lambda})$ annihilates $L(S_{\lambda}^*)\cong
L(S_{\lambda})$ there is an induced antiautomorphism on
$H_c/M_{\lambda}H_c$. Under this antiautomorphism
$\mathbb{C}[\mathfrak{h}^*]\cdot \epsilon_{\lambda}'$ is sent to
$\epsilon_{\lambda}\cdot \mathbb{C}[\mathfrak{h}]$. Since all
$\mathfrak{S}_n$--modules are self--dual, we deduce that
$$\epsilon_{\lambda}\cdot \mathbb{C}[\mathfrak{h}] \cong
\Ind_{\mathfrak{S}_{\lambda}}^{\mathfrak{S}_n}\mathbf{sign}.$$

The final sentence of the theorem is a restatement of Young's
construction of the irreducible representations of
$\mathfrak{S}_n$, \cite[Corollary 4.16]{gin}.
\end{proof}
\subsection{} We end this section by associating Theorem
\ref{3conj} to similar results found in the study of
\textit{principal nilpotent pairs}.
\subsection{} By \cite[Theorem
11.16]{EG}, there is an isomorphism of varieties
$$\spc Z_c \longrightarrow \mathcal{CM}_n,$$ where
$\mathcal{CM}_n$ denotes the Calogero--Moser space $\{ (X,Y)\in
\mat_n(\mathbb{C}) \times \mat_n(\mathbb{C}) : [X,Y] + \text{Id} =
\text{rank 1 matrix} \}//PGL_n(\mathbb{C}).$ Under this
isomorphism, the morphism $\Upsilon: \spc Z_c \longrightarrow
\mathfrak{h}^*/W \times \mathfrak{h}/W$ corresponds to the
morphism sending the matrices $(X,Y)$ to the pair of $n$--tuples
consisting of their eigenvalues. An explicit description of the
bijection between the partitions of $n$ and the elements of
$\mathcal{CM}_n$ in $\Upsilon^{-1}(0)$ is given in \cite[Section
6]{wil}.
\subsection{}
Let $\hil$ be the Hilbert scheme of $n$ points in the plane, whose
points are codimension $n$ ideals in $\mathbb{C}[X,Y]$, see
\cite{nak} for details. The Chow morphism $\hil \longrightarrow
(\mathfrak{h}^*\oplus \mathfrak{h})/W=\symn$ sends such an ideal
to its support (counted with multiplicity). There is a
$\mathbb{C}^*\times \mathbb{C}^*$--action on $\hil$ induced from
the $\mathbb{Z}^2$--grading on $\mathbb{C}[X,Y]$ where $X$
(respectively $Y$) is homogeneous of degree $(1,1)$ (respectively
$(1,-1)$). The $\mathbb{C}^*\times \mathbb{C}^*$--fixed points on
$\hil$ are in bijection with the partitions of $n$: on the ideal
corresponding to the fixed point of type $\lambda$, multiplication
by $X$ (respectively $Y$) yields a nilpotent transformation with
Jordan form of type $\lambda$ (respectively $\lambda'$). It
follows that these fixed points correspond to conjugacy classes of
principal nilpotent pairs in $\mathfrak{sl}_n(\mathbb{C})$,
\cite[Section 5]{gin}.
\subsection{}
Nakajima's quiver varieties relate Calogero--Moser space and the
Hilbert scheme. There is a geometric family $(\mathcal{CM}_n)_t$
over $\mathbb{C}$ whose generic fibre $(\mathcal{CM}_n)_t$ for
$t\neq 0$ is $\mathcal{CM}_n$ and whose degenerate point at $t=0$
is $\hil$, \cite[Section 3.2]{nak}. This family admits a
$\mathbb{C}^*\times \mathbb{C}^*$--action, which agrees with the
above $\mathbb{C}^*\times \mathbb{C}^*$--action on $\hil$ when
$t=0$. For $t\neq 0$, the diagonal $\mathbb{C}^*$--action dilates
$t$ and the anti--diagonal action corresponds to the
$\mathbb{C}^*$--action on $Z_{c}$ considered in \ref{defh}. Fixing
a partition $\lambda\vdash n$, the anti--diagonal
$\mathbb{C}^*$--fixed points corresponding to $\lambda$ give a
section of the family $(\mathcal{CM})_t$ over $\mathbb{C}\setminus
\{ 0\}$. The closure of this section yields the
$\mathbb{C}^*\times \mathbb{C}^*$--fixed point of $\hil$
corresponding to $\lambda$ and hence a representative
$(X_{\lambda}, Y_{\lambda})$ in the conjugacy class of principal
nilpotent pairs of type $\lambda$ of
$\mathfrak{sl}_n(\mathbb{C})$. The Springer correspondence
associates to $X_{\lambda}$ (respectively $Y_{\lambda}$) the
simple module $S_{\lambda}$ (respectively $S_{\lambda'}$),
\cite[Corollary 4.16]{gin}. This confirms (an amended form) of
\cite[Conjecture 17.13]{EG}.

\section{Remarks}
\subsection{} A basic problem in the representation
theory of $H_{c}$ is:

\smallskip\textbf{Problem 1:} {Find the graded $W$--character of
$L(S)$ for all $S\in \Lambda$.} \smallskip \newline  This simply
means finding the composition multiplicities in $[L(S)] =
\sum_{T\in \Lambda} (L(S): T[i])q^i[T]$. By Lemma \ref{vermmul}
this is equivalent to the problem:

 \smallskip \textbf{Problem 1$'$:} {For all $S\in \Lambda$ find the polynomials $m_{S,T}(q) =
 \sum_{T\in \Lambda}(M(S): L(T)[i])q^i$.}
\smallskip \newline Obviously, to solve Problem 1$'$ we can work in the blocks of the algebra
$\overline{H_c}$. Thus an important subproblem is to describe
which simple modules lie in a given block:

\smallskip \textbf{Problem 2:} Determine the partition of simple $W$--modules
yielding the blocks of $\overline{H_c}$.
\smallskip \newline In the notation of (\ref{the}) this corresponds to determining the fibres
of $\Theta$.

In Section 4 we have solved Problems 1$'$ and 2 for simple
$W$--modules, $S$, such that $\Theta(S)$ is a smooth point of
$\spc (Z_c)$ and in Section 5 we explictly solved Problem 1.

\subsection{} \label{brown}Problem 2 is closely related to the determination of
which complex reflection groups have $\spc (Z_c)$ smooth for
generic values of $c$. Thanks to Section 4, a point of
$M\in\Upsilon^{-1}(0)$ is singular in $\spc (Z_c)$ if
$\Theta^{-1}(M)$ is not a singleton. The following lemma, due to
Brown, proves the converse, showing that an arbitrary point $M'\in
\spc (Z_c)$ is Azumaya if and only if there is a unique simple
$H_c$-module annihilated by $M'$. One might hope that if the
points of $\Upsilon^{-1}(0)$ are all smooth then $\spc (Z_c)$ is
smooth.
\begin{lem}[K.A.Brown] Let $\mathfrak{m}$ be a maximal ideal of $H_c$.
Suppose there exists a unique maximal ideal $M$ of $H_c$ with
$M\cap Z_c = \mathfrak{m}$. Then $\mathfrak{m}$ is an Azumaya
point of $\spc Z_c$.
\end{lem}
\begin{proof}
Let $H=H_c$ and $Z=Z_c$. Since $M$ is the unique maximal ideal of
$H$ lying over $\mathfrak{m}$ we can localise at $M$ in $H$ to
obtain $H_M$. Moreover $H\subseteq H_M$ since $H$ is prime. Let
$\mathbf{e}$ be the trivial idempotent of $\mathbb{C}W$,
considered as an element of $H$. Since $H_M$ is a local ring, by
\cite[Theorem 1]{fulshut} there exists a positive integer $t$ such
that $H_M\mathbf{e} = P^{\oplus t}$, where $P$ is the unique
indecomposable projective $H_M$--module. But, by \cite[Theorem
3.1]{EG} $\mathbf{e}H\mathbf{e} \cong Z$, so that
$\mathbf{e}H_M\mathbf{e} \cong Z_\mathfrak{m}.$ Thus
$\edo_{H_M}(H_M\mathbf{e}) \cong Z_{\mathfrak{m}}$, showing that
$t=1$. Since $H_M \cong P^{\oplus n}$ for some $n$ we have
$$H_M \cong \edo_{H_M}(H_M) \cong \mat_n(Z_{\mathfrak{m}}).$$
Therefore $Z_M$ and $H_M$ are Morita equivalent, and since $H_M$
has finite global dimension so too does $Z_{\mathfrak{m}}$. It
follows that $\mathfrak{m}$ is a smooth point of $\spc Z_c$. The
lemma follows from \ref{azsm}.
\end{proof}
\subsection{} The analysis of Section 4 does provide a list of new
examples for which $\spc(Z_c)$ is not smooth at generic values of
$c$. The following results generalises the dihedral case studied
in \cite[Section 16]{EG}.  \begin{prop} Suppose $W$ is a finite
Coxeter group of type $D_{2n}$ ($n\geq 2$), $E$, $F$, $H$ or
$I_2(m)$ ($m\geq 5$). Then $Z_c$, the centre of $H_c$, is singular
for all values of $c$.
\end{prop}
\begin{proof} Suppose there exists a value $c$ such that $Z_c$ is
smooth. Thanks to Corollary \ref{maincor}(1) the blocks of
$\overline{H_c}$ are in bijective correspondence with the simple
$W$--modules and by \ref{azsm} each simple $H_c$--module, $L(S)$,
is isomorphic to the regular representation of $W$. The proof of
Theorem \ref{et} shows that $(M(S):L(S))$ equals the coefficient
of $q^{b_S}$ in the fake degree $f_{S^{\ast}}(q)$. Since Lemma
\ref{vermmul} shows that $(M(S):L(S)) = 1$ this coefficient is
necessarily $1$ for all fake degrees. But there fake degrees for
$D_{2n}, E_7$ and $E_8$ which have coefficient 2, \cite[Remark
5.6.7]{geck}, \cite{benlus}.

Let $j_S$ be the degree of $f_{S^{\ast}}(q)$, so that the trivial
module appears in $M(S)_{j_S}$. Since $L(S)$ is the only
composition factor of $M(S)$ we thus find a bound on the dimension
of $L(S)$ given by
\begin{equation} \label{bound} \dim L(S) \leq \dim(S) \left( \dim \mathbb{C}[\mathfrak{h}^*]^{co W}_{<{b_S}} + \dim
\mathbb{C}[\mathfrak{h}^*]^{co W}_{>j_S} +\min \{ \dim
\mathbb{C}[\mathfrak{h}^*]^{co W}_{b_S}, \dim
\mathbb{C}[\mathfrak{h}^*]^{co W}_{j_S}\}\right).\end{equation}
Using \cite{benlus} and \cite{alvlus} we can calculate the right
hand side of (\ref{bound}) and hence show that there exists a
simple $W$--module $S$ such that $\dim L(S)< |W|$ in types $E_6,
F_4, H_3, H_4$ and $I_2(m)$ ($m\geq 5$).

For $W$ the Coxeter group of type $E_6$ (respectively $F_4$,
$H_4$, $I_2(m)$ for $m\geq 5$) we take $S$ to be the unique
10--dimensional simple module (respectively the unique
12--dimensional simple module, the representation labelled
$\chi_7$ in \cite{alvlus}, the reflection representation); the
right hand side of (\ref{bound}) equals $49960$ (respectively
$1020$, $8808$, $8$) which is less than the order of $W$.

For $W$ the Coxeter group of type $H_3$ let $S$ be the reflection
representation, which has fake degree $f_S(q)=q+q^5+q^9$. In this
case the right hand side of (\ref{bound}) equals $120$, the order
of the group $W$. Let $S'=S\otimes \ep$ be the tensor product of
$S$ with the sign representation. The fake degree of $S'$ is
$q^{(9+5+1)}f_s(q^{-1})= q^{14}+q^{10}+q^6$. Thus there are two
copies of the sign representation in $M(S)_{>9}$, one in degree
10, the other in degree 14. Hence we can subtract $1$ from our
estimate in (\ref{bound}), showing that $\dim L(S)< 120=|W|$, as
required.
\end{proof}
There are a number of other complex reflection groups whose fake
degrees have leading coefficient greater than $1$, so we can prove
non--smoothness in greater generality. The techniques above do not
work for $D_{2n+1}$ however.

\subsection{} Type $G_2$ is the smallest example of a Weyl group
for which $Z_c$ is singular for generic values of $c$. In this
case, calculation shows that $\Upsilon^{-1}(0)$ has five elements,
with $D_{\pm}$, the pair of two dimensional simple modules, being
the only simple $W$--modules sent to the same element under
$\Theta$. The polynomials of Problem 1$'$ are $$ m_{D_+,D_+}(q) =
m_{D_-,D_-}(q) = 1 + q^4, \quad m_{D_+,D_-}(q)=m_{D_-,D_+}(q) =
q+q^3.$$

\end{document}